%%%%%%%%%%%%%%%%%%%%%%%%%%%%%%%%%%%%%%%%%%%%%%%%%%%%
% typoref.tex. V : January 18, 2000.
% Author : Anthony PHAN
% Warning : syntaxe +- LaTeX
% Sources :
% T. Lachand--Robert, ``La Ma\^\i trise de \TeX'',
% R\'ef\'erences crois\'ees;
% latex.ltx's sources;
% and of course the \TeX book.
%%%%%%%%%%%%%%%%%%%%%%%%%%%%%%%%%%%%%%%%%%%%%%%%%%%%%
%

\catcode`\@=11

\magnification=1200
\baselineskip=14pt

\pretolerance=500    \tolerance=1000 \brokenpenalty=5000

\catcode`\;=\active
\def;{\relax\ifhmode\ifdim\lastskip>\z@
\unskip\fi\kern.2em\fi\string;}

\overfullrule=0mm

\catcode`\!=\active
\def!{\relax\ifhmode\ifdim\lastskip>\z@
\unskip\fi\kern.2em\fi\string!}

\catcode`\?=\active
\def?{\relax\ifhmode\ifdim\lastskip>\z@
\unskip\fi\kern.2em\fi\string?}

\frenchspacing

\newif\ifpagetitre            \pagetitretrue
\newtoks\hautpagetitre        \hautpagetitre={ }
\newtoks\baspagetitre         \baspagetitre={1}

\newtoks\auteurcourant        \auteurcourant={ }
\newtoks\titrecourant
\titrecourant={ }
\newtoks\hautpagegauche       \newtoks\hautpagedroite
\hautpagegauche={\hfill\sevenrm\the\auteurcourant\hfill}
\hautpagedroite={\hfill\sevenrm\the\titrecourant\hfill}

\newtoks\baspagegauche       \baspagegauche={\hfill\rm\folio\hfill}

\newtoks\baspagedroite       \baspagedroite={\hfill\rm\folio\hfill}

\headline={
\ifpagetitre\the\hautpagetitre
\global\pagetitrefalse
\else\ifodd\pageno\the\hautpagedroite
\else\the\hautpagegauche\fi\fi}

\footline={\ifpagetitre\the\baspagetitre
\global\pagetitrefalse
\else\ifodd\pageno\the\baspagedroite
\else\the\baspagegauche\fi\fi}

\def\date{\ {\the\day}\
\ifcase\month\or Janvier\or F\'evrier\or Mars\or Avril
\or Mai \or Juin\or Juillet\or Ao\^ut\or Septembre
\or Octobre\or Novembre\or D\'ecembre\fi\
{\the\year}}

\def\up#1{\raise 1ex\hbox{\sevenrm#1}}

\def\cqfd{\unskip\kern 6pt\penalty 500
\raise -2pt\hbox{\vrule\vbox to 10pt{\hrule width 4pt
\vfill\hrule}\vrule}\par\medskip}

\def\section#1{\vskip 7mm plus 20mm minus 1.5mm\penalty-50
\vskip 0mm plus -20mm minus 1.5mm\penalty-50
{\bf\noindent#1}\nobreak\smallskip}

\def\subsection#1{\medskip{\bf#1}\nobreak\smallskip}

\def\displaylinesno #1{\dspl@y\halign{
\hbox to\displaywidth{$\@lign\hfil\displaystyle##\hfil$}&
\llap{$##$}\crcr#1\crcr}}

\def\ldisplaylinesno #1{\dspl@y\halign{
\hbox to\displaywidth{$\@lign\hfil\displaystyle##\hfil$}&
\kern-\displaywidth\rlap{$##$}
\tabskip\displaywidth\crcr#1\crcr}}

\def\hfl#1#2{\smash{\mathop{\hbox to 12 mm{\rightarrowfill}}
\limits^{\scriptstyle#1}_{\scriptstyle#2}}}

%
% style (look at the behavior of \item dans \bibitem too,
% and at one ,\  in \re@dreferenceslist)
% Feel free to change: 	\bibn@me (title like ``R\'ef\'erences'')
%			\bibliographym@rk (general style)
%
\def\bibn@me{R\'ef\'erences}
\def\bibliographym@rk{\centerline{{\sc\bibn@me}}
	\sectionmark\section{\ignorespaces}{\unskip\bibn@me}
	\bigbreak\bgroup
	\ifx\ninepoint\undefined\relax\else\ninepoint\fi}
%
% Beware of the \bgroup: it will be closed by \endthebibliography
%
% \refsp@ce is the spacing command that appens between multiple
% references.
%
\let\refsp@ce=\
\let\bibleftm@rk=[
\let\bibrightm@rk=]
%
% if you want more space between brackets...
%\let\refsp@ce=\thinspace
%\def\bibleftm@rk{[\thinspace}
%\def\bibrightm@rk{\thinspace]}
%
% frenchy stuff
%
\def\numero{n\raise.82ex\hbox{$\fam0\scriptscriptstyle
o$}~\ignorespaces}
%
% new variables
%
\newcount\equationc@unt
\newcount\bibc@unt
\newif\ifref@changes\ref@changesfalse
\newif\ifpageref@changes\ref@changesfalse
\newif\ifbib@changes\bib@changesfalse
\newif\ifref@undefined\ref@undefinedfalse
\newif\ifpageref@undefined\ref@undefinedfalse
\newif\ifbib@undefined\bib@undefinedfalse
\newwrite\@auxout
%
% mark an equation
%
%\def\eqnum{\global\advance\equationc@unt by 1%
%\edef\lastref{\number\equationc@unt}%
%\eqno{(\lastref)}}
%
% One can reference anything, just copy the former macro
% and use it so: \machin \label{truc}
% In machin you would have defined \lastref by some number
% or any text.
%
% References macros
%
% The next macros are the core of \ref and \cite commands.
% Its first argument may be ref, pageref or bib.
%
% It is too tricky to be explained.
% (It is a bit recursive.)
% It allows using \cite or \ref or ...
% with arbitrary many arguments,
% for instance:
% \cite{knuth1,knuth2,ma pomme}
%
% First argument is always ref, pageref or bib.
%
\def\re@dreferences#1#2{{%
	\re@dreferenceslist{#1}#2,\undefined\@@}}
\def\re@dreferenceslist#1#2,#3\@@{\def\next{#2}%
	\expandafter\ifx\csname#1@@\meaning\next\endcsname\relax
	??\immediate\write16
	{Warning, #1-reference "\next" on page \the\pageno\space
	is undefined.}%
	\global\csname#1@undefinedtrue\endcsname
	\else\csname#1@@\meaning\next\endcsname\fi
	\ifx#3\undefined\relax
	\else,\refsp@ce\re@dreferenceslist{#1}#3\@@\fi}
%
% notice that the former ``,\refsp@ce'' will separate
% multiple arguments. But beware of spaces
% while defining a reference or calling for it!
%
% tricky thing: \newlabel has two arguments
% {labelname}{{\lastref}{\pageref}}
% The second argument is read as two arguments
% by \newl@bel. This was necessary to get
% a jobname.aux containing the same syntax
% LaTeX would produce and use.
%
\def\newlabel#1#2{{\def\next{#1}\newl@bel#2}}
\def\newl@bel#1#2{%
	\expandafter\xdef\csname ref@@\meaning\next\endcsname{#1}%
	\expandafter\xdef\csname pageref@@\meaning\next\endcsname{#2}}
\def\label#1{{%
	\toks0={#1}\message{ref(\lastref) \the\toks0,}%
	\ignorespaces\immediate\write\@auxout%
	{\noexpand\newlabel{\the\toks0}{{\lastref}{\the\pageno}}}%
	\def\next{#1}%
	\expandafter\ifx\csname ref@@\meaning\next\endcsname\lastref%
	\else\global\ref@changestrue\fi%
	\newlabel{#1}{{\lastref}{\the\pageno}}}}
\def\ref#1{\re@dreferences{ref}{#1}}
\def\pageref#1{\re@dreferences{pageref}{#1}}
%
% bibliography macros
%
\def\bibcite#1#2{{\def\next{#1}%
	\expandafter\xdef\csname bib@@\meaning\next\endcsname{#2}}}
\def\cite#1{\bibleftm@rk\re@dreferences{bib}{#1}\bibrightm@rk}
%
% The argument of \beginthebibliography
% is any sequence of numerals which will represent
% the maximum \item's length. If you have less than 9
% \bibitem's, this argument may be {any numeral}.
% if you have between 100 and 999 \bibitem's
% this argument may be {any three numerals},
% and so on.
%
\def\beginthebibliography#1{\bibliographym@rk
	\setbox0\hbox{\bibleftm@rk#1\bibrightm@rk\enspace}
	\parindent=\wd0
	\global\bibc@unt=0
	\def\bibitem##1{\global\advance\bibc@unt by 1
		\edef\lastref{\number\bibc@unt}
		{\toks0={##1}
		\message{bib[\lastref] \the\toks0,}%
		\immediate\write\@auxout
		{\noexpand\bibcite{\the\toks0}{\lastref}}}
		\def\next{##1}%
		\expandafter\ifx
		\csname bib@@\meaning\next\endcsname\lastref
		\else\global\bib@changestrue\fi%
		\bibcite{##1}{\lastref}
		\medbreak
		\item{\hfill\bibleftm@rk\lastref\bibrightm@rk}%
		}
	}
\def\endthebibliography{\egroup\par}
%
% THE NEXT MACRO MUST BE INCLUDED
% IN THE \BYE COMMAND. FOR INSTANCE:
%
   %\catcode`@=11
   \outer\def\bye{\@closeaux
   	\par\vfill\supereject\end}
   %\catcode`@=12
%
\def\@closeaux{\closeout\@auxout
	\ifref@changes\immediate\write16
	{Warning, changes in references.}\fi
	\ifpageref@changes\immediate\write16
	{Warning, changes in page references.}\fi
	\ifbib@changes\immediate\write16
	{Warning, changes in bibliography.}\fi
	\ifref@undefined\immediate\write16
	{Warning, references undefined.}\fi
	\ifpageref@undefined\immediate\write16
	{Warning, page references undefined.}\fi
	\ifbib@undefined\immediate\write16
	{Warning, citations undefined.}\fi}
%
% initialization of jobname.aux
%
\immediate\openin\@auxout=\jobname.aux
\ifeof\@auxout \immediate\write16
    {Creating file \jobname.aux}
\immediate\closein\@auxout
\immediate\openout\@auxout=\jobname.aux
\immediate\write\@auxout {\relax}%
\immediate\closeout\@auxout
\else\immediate\closein\@auxout\fi
%
% Let's read this file and open it out
%
\input\jobname.aux \par
\immediate\openout\@auxout=\jobname.aux
% this file will be closed by \bye.
%
% That's all, folks!
%

\def\Z{{\bf Z}} 
\def\R{{\bf R}}

\def\bP{{\bf P}}
\def\bQ{{\bf Q}}
\def\bR{{\bf R}}

  \def\pro{\noindent {\bf{Proof :
}}}
   \def\resp{{\rm resp.  }}

   \def\la{{\lambda}}

\def\om{{\omega}}
\def\omc{{\hat{\omega}}}
\def\and{\quad\hbox{and}\quad}

 \def\ux{{\underline{x}}}

 \def\ens{\enspace} \def\noi{\noindent}

\def\build#1_#2^#3{\mathrel{\mathop{\kern 0pt#1}\limits_{#2}^{#3}}}

%\newfam\gothfam \scriptscriptfont\gothfam=\fivegoth
%\textfont\gothfam=\tengoth \scriptfont\gothfam=\sevengoth
%\def\goth{\fam\gothfam\tengoth}

%\def\cqfd{\unskip\kern 6pt\penalty 500 \raise 0pt\hbox{\vrule\vbox
%to6pt{\hrule width 6pt \vfill\hrule}\vrule}\par}

\def\pro{\noindent {\bf Proof : }}

\def\smallsquare{\vbox{\hrule\hbox{\vrule height 1 ex\kern 1
ex\vrule}\hrule}}
\def\cqfd{\hfill \smallsquare\vskip 3mm}

\def\hw{{\hat w}}
\def\hla{{\hat \lambda}}
\def\utheta{{\underline{\theta}}}

\def\bibn@me{R\'ef\'erences bibliographiques}
%\input typpo
%
%\catcode`@=11
\def\bibliographym@rk{\bgroup}
%
% \bye est modifie pour la biblio et la table des matieres
%
\outer\def\bye{ 	\par\vfill\supereject\end}

%%%%%%%%%%%%%%%%%%%%%%%%%%%%%%%%%%%%%%%%%%
%\catcode`@=12

\null

%\centerline{}

\vskip 2mm

\centerline{\bf Exponents of Diophantine approximation}

\vskip 8mm

\centerline{Yann B{\sevenrm UGEAUD} \footnote{}{\rm
}
\& Michel L{\sevenrm AURENT}
}

\vskip 11mm

\section{ 1. Introduction}

The well known Dirichlet Theorem asserts that for any irrational
real number $\xi$ and {\it any} real number $Q \ge 1$, there exist
integers $p$ and $q$ with $1 \le q \le Q$ and
$$
|q \xi - p| \le Q^{-1}.  \eqno (1)
$$
As observed by Khintchine \cite{Kh26b}, there is no $\xi$ for
which the exponent of $Q$ in (1) can be lowered (see
Davenport \& Schmidt \cite{DSa} or Schmidt \cite{SchmLN}
for a very precise result).
However, for any $w>1$, there clearly exist real numbers $\xi$
for which, for {\it arbitrarily large} values of $Q$, the equation
$$
|q \xi - p| \le Q^{-w}
$$
has a solution in integers $p$ and $q$ with $1 \le q \le Q$.
Obviously, the quality of approximation strongly depends on
whether we are interested in a uniform statement (i.e., a statement
valid
for any $Q$, or for any $Q$ sufficiently large) or in a statement
valid for arbitrarily large $Q$.
In the case of rational approximation, these
questions are quite well understood, essentially thanks to the
continued fraction theory.
However, the Dirichlet Theorem extends well to rational simultaneous
approximation, and to simultaneous approximation of linear forms.
In Section 2, we define exponents of Diophantine approximation
related to these questions, and survey known results on them.
The two dimensional case is now fully understood and the results are
displayed in Section 3. In order to extend possibly these results in
higher  dimension,
we define geometrically  in Section 4 more  exponents of Diophantine
approximation,
in connection with the work of Schmidt \cite{SchmB}.

It is a notorious fact that questions of Diophantine approximation are
in general much more difficult when the quantities we
approximate are dependent. Classical examples include the simultaneous
rational approximation of the first $n$ powers of a transcendental
number,
and the approximation of linear forms whose coefficients are precisely
the first $n$ powers of a transcendental number.
These are considered in Section 5.
When $n=2$, important progress has been
made recently by Roy \cite{RoyA,RoyB,RoyC,RoyD,RoyE,RoyF},
and Section 6  is devoted to his
results and some of their extensions.

The present paper is in great part a survey, however, it contains
several new results (e.g., Theorems 8 and 10).
General references for the topic investigated here
are the seminal paper of Khintchine \cite{Kh26b}, Cassels' book
\cite{Cas}
and the monograph \cite{BuLiv}.

\section{ 2. Approximation of independent quantities}

We begin with some notations and definitions.
If  $\utheta$ is a  (column) vector  in $\R^n$, we denote by $\vert
\utheta\vert$
the  maximum of the   absolute values of its coordinates  and by
$$
\Vert \utheta\Vert = \min_{\ux \in \Z^n} \vert \utheta -\ux\vert
$$
the maximum of the distances of its  coordinates to the  rational
integers.

Following the convention introduced in \cite{BuLaA}, we indicate by
a `hat' the exponents of {\it uniform} Diophantine approximation.

\proclaim Definition 1. Let $n$ and $m$ be positive integers
and let $A$ be a real matrix with $n$ rows and $m$ columns.
We denote by  $\om_{n, m}(A)$
the supremum of the real numbers $w$ for which, 
{\it for arbitrarily large} 
real numbers $X$, the inequalities
$$
\Vert A\ux \Vert \le X^{-w}
\quad {\it and} \quad \vert \ux \vert \le X   \eqno{(2)}
$$
have a non-zero solution  $\ux$ in $\Z^m$. We
denote by  $\omc_{n, m}(A)$ the supremum of the real numbers $w$
for which, {\it for all sufficiently large} positive real numbers  $X$,  
the inequalities (2) have a non-zero integer solution $\ux$ in $\Z^m$.

For a $n \times m$ matrix $A$, the Dirichlet box principle implies that
$$
\om_{n, m} (A) \ge  \omc_{n, m}(A) \ge {m\over n}.   \eqno (3)
$$
Furthermore, we have both equalities in (3) for almost all
matrices $A$, with respect to the Lebesgue measure on $\R^{mn}$, as
follows from the Borel--Cantelli Lemma.

The left-hand side inequality of (3) has been improved by
Jarn\'\i k \cite{Jar50,Jar54} as follows.

\proclaim Theorem 1. For any $n \ge 2$ and any $n\times 1$ real matrix
$A$ with at least two coefficients which are $\bQ$-linearly 
independent modulo $\bQ$, we have
$$
\om_{n,1} (A) \ge {\omc_{n, 1}^2 (A) \over 1 - \omc_{n, 1} (A)}.
\eqno{(4)}
$$
For any $n \ge 1$ and any $n\times 2$ real matrix $A$, we have
$$
\om_{n,2} (A) \ge \omc_{n,2} (A) (\omc_{n,2} (A) - 1). \eqno(5)
$$
For any $n \ge 1$, any $m \ge 3$ and any $n\times m$ real matrix $A$
with $\omc_{n,m} (A) > (5m^2)^{m-1}$, we have
$$
\om_{n,m} (A) \ge \bigl(\omc_{n,m} (A) \bigr)^{m/(m-1)}
- 3 \, \omc_{n,m} (A).
$$

In all what follows, we denote by ${}^t A$ the transpose of
the matrix $A$.
It is well-known that $\om_{n,m}(A)$ and $\om_{m,n} ({}^t A)$
are linked by a transference principle. Dyson \cite{Dy47} established
the lower bound
$$
\om_{n,m}(A) \ge {m \, \om_{m,n}({}^t A) + m - 1 \over
(n - 1) \om_{m,n}({}^t A) + n}, \eqno (6)
$$
thus extending earlier results of
Khintchine \cite{Kh25,Kh26b} who delt with the case $\min\{n ,m\} = 1$.
For a proof, the reader is referred to
Gruber \& Lekkerkerker \cite{GrLe}, Section 45.3,
Cassels \cite{Cas}, Chapter V, Theorem IV, or Schmidt
\cite{SchmLN}, Chapter IV, Section 5.
Inequalities (6) have been shown to be best possible for $\min\{n ,m\}
= 1$
by Jarn\'\i k \cite{Jar35,Jar36}, who also got some related results
\cite{Jar35a}. For general $m$ and $n$, Jarn\'\i k \cite{Jar59}
proved that (6) is best possible except, possibly, when
$1 < n < m$ and $\om_{n,m} (A) < (m-1)/(n-1)$, in which case his method
does not give anything.

Furthermore, extending earlier results of
Jarn\'\i k \cite{Jar38}, Apfelbeck \cite{Apf} established that the
uniform
exponents $\omc_{n,m} (A)$ and  $\omc_{m,n} ({}^t A)$
are linked by the same relation
$$
\omc_{n,m} (A) \ge {m \, \omc_{m,n} ({}^t A) + m - 1
\over (n-1) \omc_{m,n} ({}^t A) + n}. \eqno (7)
$$
Jarn\'\i k \cite{Jar38} and Apfelbeck \cite{Apf} succeeded in
improving (7) when either $\omc_{n,m} (A)$ or $\omc_{m,n} ({}^t A)$
is large.
Before investigating the set of values taken by
the functions ${\om}_{n, m}$ and $\omc_{n, m}$, we introduce
the following definition.

\proclaim Definition 2.
By {\it spectrum} of an exponent of Diophantine approximation, we mean
the
set of values taken by this exponent on the set of
$n$ by $m$ real matrices $A$ of maximal rank, i.e. of rank $\min\{m,
n\}$.

Since $\om_{1,1}((\xi)) = \om_{2,1} ({}^t (\xi, \xi))$
holds for any real number $\xi$, we give the above definition of
spectrum
to avoid trivialities.

Except for $m=n=1$ (in that case, we can use the continued fraction
theory),
it is in general a difficult problem
to construct explicit examples of regular $n$ by $m$ matrices $A$
with prescribed values for ${\om}_{n, m}(A)$
and/or for $\omc_{n, m}(A)$. However, the spectrum of the function
$\om_{n, m}$ has been completely determined, thanks to a deep result
of Dickinson \& Velani \cite{DV} (see also \cite{BDV}), who calculated
the
Hausdorff dimension of the set of matrices $A$ with ${\om}_{n,
m}(A)=\tau$,
for an arbitrary real number $\tau$.

\proclaim Theorem 2. For any positive integers $n$ and $m$, the 
spectrum
of the function $\om_{n, m}$ is equal to $[m/n, + \infty]$.

As for the spectra of the exponents $\omc_{n, m}$, much less is known.
They are contained in $[1/n, 1]$ if $m=1$ and in
$[m/n, + \infty]$ if $m\ge 2$ (this is an immediate consequence of 
(3)).
In particular, we have $\omc_{1, 1} ((\xi)) = 1$ for any
irrational real number $\xi$.
The situation is completely different in the case
$(m, n) \not= (1, 1)$.

\proclaim Theorem 3.
For any positive integers $m$, $n$ with $m \ge 2$
there are continuum many $n$ by $m$ real matrices $A$
whose coefficients are algebraically independent
and which satisfy $\omc_{n, m}(A) = + \infty$.
For any positive integer $n$, there are continuum many
$n$ by $1$ real matrices $A$
whose coefficients are algebraically independent
and which satisfy $\omc_{n, 1}(A) = 1$.

\pro For $(m, n) = (2, 1)$ or $(1, 2)$ and the coefficients of the
corresponding matrices are linearly independent, this is
due to Khintchine \cite{Kh26b}
(see also Theorem XIV, page 94, of \cite{Cas}).
Further results have been obtained by Chabauty \& Lutz \cite{ChLu}.
Jarn\'\i k \cite{Jar59b} completed the proof of the theorem, using
a quite different approach (see also Lesca \cite{Les}). \cqfd

We address the following problem, which is likely to be difficult.

\proclaim Problem 1. For positive integers $m$ and $n$,
determine the spectrum of the function $\omc_{n, m}$.

Partial results when $\min \{n, m\} = 1$ have been established
by Jarn\'\i k \cite{Jar54}.

\proclaim Theorem 4. The spectrum of $\omc_{1,2}$ (resp. $\omc_{2,1}$)
is equal to $[2, +\infty]$ (resp. to $[1/2,1]$). For any integers
$m \ge 2$ and $n\ge 2$, the spectrum of $\omc_{1,m}$
contains $(2^{m-1}, + \infty]$ and that of $\omc_{n,1}$
contains $((u_n-2 - u_n^{-n+1})/(u_n-1),1]$, where $u_n$ is the largest
real
root of the polynomial $X^{n-1} - X^{n-2} - \sum_{k=0}^{n-2} \, X^k$.

Jarn\'\i k's proof of Theorem 4 is constructive and rests on the
continued fraction theory.

\medskip

In the above definition of the exponent
$\om_{n, m}$ (\resp ${\omc}_{n, m}$),
we do not require that there exists a positive constant $c$
such that, for arbitrarily large
real numbers $X$ (\resp for any sufficiently large real number $X$),
the inequalities
$$
\Vert A\ux \Vert \le c \, X^{-\om_{n, m}}
\quad {\rm and} \quad \vert \ux \vert \le X
$$
and
$$
\Vert A\ux \Vert \le c \, X^{-\omc_{n, m}}
\quad {\rm and} \quad \vert \ux \vert \le X,
$$
respectively, have a non-zero solution $\ux$ in $\Z^m$.
Taking this into consideration yields new problems.

Actually, Dirichlet's Theorem implies that, for any $X > 1$,
the inequations %%m
$$
\Vert A\ux \Vert \le c \, X^{-m/n}
\quad {\rm and} \quad \vert \ux \vert \le X   \eqno (8)
$$
have a non-zero solution  $\ux$ in $\Z^m$, when $c=1$. This suggests to
us
to introduce the following definition.

\proclaim Definition 3. Let $A$ be a $n \times m$ real matrix. We say
that
Dirichlet's Theorem can be improved for the matrix $A$ if there exists
a positive constant $c < 1$ such that (8) has a solution  $\ux$ in
$\Z^m$
for any sufficiently large $X$.

When $m=n=1$, that is, when $A = ((\xi))$ for some irrational
real number $\xi$, it is well known that Dirichlet's Theorem can be
improved if, and only if, $\xi$ has bounded partial quotients in its
continued fraction expansion. A precise statement has been
obtained by Davenport \& Schmidt \cite{DSa}. In particular, the set of
$1\times 1$ matrices $A$ for which Dirichlet's Theorem can be %%m
improved has Lebesgue measure zero and Hausdorff dimension 1. %%m
This assertion  has been partly extended  %%m
to linear forms and to simultaneous approximation
by Davenport \& Schmidt \cite{DaScB}.

\proclaim Theorem 5. For any positive integer $n$, the set of
$n\times 1$ (resp.  of $1\times n$)  matrices for which Dirichlet's
Theorem
can be improved has $n$-dimensional Lebesgue measure zero.

Notice that Khintchine \cite{KhB} %%m
proved  that the set of {\it singular} $n\times m$ real matrices $A$ (meaning that    %%m
for {\it all }  positive constant $c$, the inequations (8) have a  non-zero solution $\ux$ in $\Z^m$ %%m
for all   $X$ greater than $X_0(A,c)$),   has %%m
$mn$-dimensional Lebesgue measure zero. This weaker result is a  consequence of the %%m
Borel--Cantelli lemma (see \cite{Cas}, page 92). %%m
The proof of Theorem 5 is quite involved. %%m
According to Kleinbock and Weiss \cite{KlWe}, it can be  %%m
generalized to $n \times m$ matrices. Actually, a more general  %%m
result is proved in \cite{KlWe}.    Maybe, it is possible to  adapt the methods of \cite{DaScB,KlWe}  %%m
to solve the following problem, which seems to be rather difficult.  %%m

\proclaim Problem 2. Let $c$ be a real number with $0<c<1$.  
Determine the Hausdorff dimension    %%m
of the set of $n\times m$ matrices such that (8) has a solution  $\ux$ in  %%m
$\Z^m$ for any sufficiently large $X$.

We end this section by briefly mentionning that we can as well
consider inhomogeneous problems in Diophantine approximation (see
Chapters III and V from \cite{Cas}). The corresponding
exponents of approximation have been introduced in \cite{BuLaB}, where
it is
established that the
exponent of approximation to a generic point in $\R^n$ by a system of
$n$ linear forms is equal to the inverse of the
uniform homogeneous exponent associated
to the system of dual linear forms.

\section{ 3. Approximation in dimension two}

We investigate  more precisely in this section the above spectra when
$A$ is an  $1\times 2$ or a
$2\times 1$ real matrix.
In this case, the uniform exponents $\omc(A)$ and $\omc({}^tA)$ are
linked by an equation
   due to Jarn\'\i k \cite{Jar38}.  This fact seems to
have been completely forgotten since 1938. 
%%It  has been   
%%recently rediscovered by Roy.  
%% Damien ne souhaitait pas que ceci apparaisse, je m'en souviens maintenant

\proclaim Theorem 6. For any $1\times 2$ real matrix $A=(\alpha,
\beta)$, with $\alpha$ or $\beta$ irrational,  the equality
$$
\omc_{2, 1}({}^t A) = 1 - {1 \over \omc_{1, 2} (A)}
$$
holds.

On the other hand, Khintchine's transference inequalities $(6)$ read  
here
$$
{ \om_{1,2}(A)\over \om_{1,2}(A)+2} \le \om_{2,1}({}^tA) \le { 
\om_{1,2}(A)-1\over 2},
$$
for any matrix $A$ as in Theorem 6.
Next Theorem, which is the main result of \cite{LauB}, refines this
latter estimate.

\proclaim
Theorem 7.
For any row vector  $A = (\alpha, \beta) $ with $1,\alpha , \beta$
linearly independent over $\bQ$,  the four exponents
$$
v= \om_{1,2}(A),\quad  v'= \om_{2,1}({}^t A),
\quad  w = \omc_{1,2}(A), \quad w'= \omc_{2,1}({}^t A),
$$
satisfy  the relations
$$
2\le w\le +\infty ,\quad  w' ={w-1\over w} , \quad {v(w-1)\over v+w}
\le v' \le {v-w+1\over w}.
$$
When  $v=+\infty$ we have to understand these relations as   $w-1 \le
v' \le +\infty$,  and $w'=1,v'=+\infty$ if
moreover  $w=+\infty$.
Conversely,  for each quadruple  $(v,v',w,w')$ in
$(\bR_{>0}\cup\{+\infty\})^4$
satisfying the previous conditions,
there exists a row vector $A=(\alpha,\beta)$ of real numbers with
$1,\alpha , \beta$ linearly independent over $\bQ$,
such that
$$
v= \om_{1,2}(A),\quad  v'= \om_{2,1}({}^t A),
\quad  w = \omc_{1,2}(A), \quad w'= \omc_{2,1}({}^t A).
$$

Notice that our  estimate
$$
{v(w-1)\over v+w} \le v' \le {v-w+1\over w}  \eqno (9)
$$
refines
Khintchine's  inequalities since $w\ge 2$.
Theorem 7 implies that the lower bound
$(5)$ is optimal when $n=1$.

\proclaim
Corollary 1.
For any row vector  $A =(\alpha, \beta)$ with $1,\alpha,\beta$ linearly
independent over $\bQ$,
the lower bounds
$$
   \omc_{1,2}(A) \ge 2 \and \om_{1,2}(A) \ge
\omc_{1,2}(A)(\omc_{1,2}(A)-1)
$$
hold. Moreover, let $v$ and $w$ be positive real numbers with
$$
w \ge 2 \and v\ge w(w-1).
$$
Then there exists a row vector $A \in \bR^2$ such that
$$
\om_{1,2}(A)=v \and \omc_{1,2}(A) =w.
$$

Similarly, the  lower bound $(4)$ of Theorem 1 is best possible for
$n=2$.

\proclaim
Corollary 2.
For any  column vector $A =\left(\matrix{\alpha\cr\beta\cr}\right)$
   with $1,\alpha,\beta$ linearly independent over $\bQ$, the
inequalities
   $$
{1\over 2} \le \omc_{2,1}(A) \le 1 \and \om_{2,1}(A) \ge
{\omc_{2,1}(A)^2\over 1 - \omc_{2,1}(A)}
$$
hold. Moreover, let $v'$ and $w' $ be positive real numbers satisfying
$$
{1\over 2 }\le w'\le 1 \and v'  \ge  { w'^2\over 1-w'}.
$$
Then there exists a column vector $A \in \bR^2$ such that
$$
\om_{2,1}(A)=v' \and \omc_{2,1}(A) =w'.
$$

For the deduction  of the corollaries, observe that, for given positive
real numbers
$v$ and $w$, the interval
$$
{v(w-1)\over v+w} \le v' \le {v-w+1\over w}
$$
occurring in Theorem 7, is non-empty exactly when $v\ge w(w-1)$.
For the minimal value $v=w(w-1)$, it reduces to the
point
$$
{(w-1)^2\over w}= {w'^2\over 1-w'}.
$$
   Corollaries 1 and 2 immediately follow.

\section{ 4. Further problems}

Having regard to Section 3, we are led to ask for an extension of
Theorem 7  in higher dimension.

\proclaim
Problem 3. Let $n$ and $m$ be positive integers with $n \le m$.
Describe the set of quadruples
$$
(\om_{n,m}(A),\om_{m,n}({}^t A),
\omc_{n,m}(A),\omc_{m,n}({}^tA)),
$$
when $A$ ranges over the set of real $n\times m$ matrices.

Even a conjectural answer to Problem 3 is unclear to us,
unless when $(m,n)=(1,1)$ or $(m,n) =(2,1)$.
A (small) contribution towards the resolution of Problem 3
is the following refinement of Khintchine's transference
inequalities for row matrices.

\proclaim Theorem 8.
Let $m \ge 2$ and $A$ be a $1 \times m$ real matrix.
Set
$$
v= \om_{1,m}(A),\quad  v'= \om_{m,1}({}^t A),
\quad  w = \omc_{1,m}(A), \quad w'= \omc_{m,1}({}^t A),
$$
Then, we have
$$
v' \ge {v \over
{(m-1) w \over w - 1} (v+1) - v} \quad {\it and} \quad
v \ge {(m-1) (1 + v') \over 1 - w'} - 1.  \eqno (10)
$$

When $m=2$, the combination of Theorems 6 and 8
yields the inequalities (9). With the notation of Theorem 8,
inequality (6) reads
$$
v' \ge {v \over
m (v+1) - v} \quad {\rm and} \quad
v \ge m (1 + v') - 1.
$$
This is weaker than (10)
since $w \ge m$ and $w' \ge 1/m$.

To establish Theorem 8, we insert in the proof of Theorem II from
Chapter V of \cite{Cas} an upper bound for the second of the
successive minima of the convex body involved. This upper bound is
obtained by following Section 9 of the same Chapter and making a
suitable use of the exponents of uniform approximation. Full
details will be given in a subsequent work.

   In order to enlighten Problem 3,
   it might be relevant to introduce  more  exponents of approximation
in intermediate dimensions. In this section, we reformulate and extend
geometrically the definitions of the exponents $\om_{n,m}(A)$ and
$\omc_{n,m}(A)$.

Let $\Theta$ and $L$ be two real linear subvarieties contained in a
projective space $\bP^N(\bR)$. Assume that $\Theta$ and $L$ are
non-empty and
distinct from $\bP^N(\bR)$, and set
$$
\delta = \dim \Theta  \and d = \dim L.
$$
Following Schmidt \cite{SchmB} (see also the appendix of \cite{Rand})
and  using  Euclidean
geometry,   we can attach to $\Theta$ and
$L$  a sequence of  sines
$$
0 \le \psi_1(\Theta,L) \le \dots \le  \psi_t(\Theta,L) \le 1, \quad
{\rm where} \quad t= \min(d+1,N-d,\delta+1,N-\delta),
$$
of various acute angles measuring in the ambient space $\bP^N(\bR)$
the proximity  of the linear subvarieties $\Theta$ and $L$. When 
$\delta +d \le N-1$,
the intersection $\Theta\cap L$ is usually empty, and the
smallest sine $\psi_1(\Theta , L)$
can be compared  with the minimal projective distance (the
normalization does not matter for our purpose)
between the  points  of $\Theta$ and $L$.
Notice that we have $t=1$
exactly when either $\Theta$ or $L$ is either a point or an hyperplane
of $\bP^{N}(\bR)$.

When $L$ is rational over $\bQ$, we denote by $H(L)$ its height, as
defined in \cite{SchmB}.  We are now able to extend Definition 1
in the following way.

\proclaim Definition 4. Let $\Theta$  be a proper  real linear
subvariety of $\bP^{N}(\bR)$ of dimension $\delta$. Let $d$ and $k$ be
integers
with
$$
0\le d \le N-1 \and 1 \le k \le \min(\delta +1,N-\delta ,d+1,N-d).
$$
We denote by $w_{d,k}(\Theta)$ the supremum of the real numbers $w$
such that for arbitrarily large real numbers $X$ there
exists a rational linear subvariety $L \subset \bP^{N}(\bR)$ of
dimension $d$, satisfying the inequations
$$
\psi_k(\Theta,L) \le H(L)^{-1}X^{-w} \and H(L)\le X. \eqno{(11)}
$$
Similarly, we denote by $\hat{w}_{d,k}(\Theta)$ the supremum of the
real numbers $w$ such that, for all sufficiently large positive real
number $X$,
there exists   a rational linear subvariety $L \subset \bP^{N}(\bR)$ of
dimension $d$ satisfying $(11)$.

The link with Definition 1 is achieved in the following way. Let $A$ be
a real $n\times m$ matrix. Put $N= m+n-1$ and associate to $A$
the $(n-1)$-dimensional linear subvariety
$$
\Theta = \bP( {\rm Span} (A\times I_n)) \subset \bP^N(\bR) =
\bP(\bR^{m+n}),
$$
whose points
have homogeneous coordinates belonging to the vectorial subspace of
$\bR^{m+n}$ spanned by the rows
of the $n\times (m+n)$ matrix $(A,I_n)$. Notice that the  application
$A\mapsto \Theta$
clearly defines a homeomorphism from $\bR^{mn}$  to an  open set of the
Grassmanian of the $(n-1)$-dimensional linear spaces
in $\bP^{N}(\bR)$.
When $d=0$, the single exponent $w_{0,1}(\Theta)$ (resp.
$\hat{w}_{0,1}(\Theta))$ measures the approximation (resp. uniform
approximation)
to $\Theta$ by rational points in $\bP^{N}(\bR)$. It is readily
observed that Definitions 1 and 4 are then consistent,  in the
sense that
$$
w_{0,1}(\Theta) = \om_{m,n}({}^tA) \and \hat{w}_{0,1}(\Theta) =
\omc_{m,n}({}^tA).
$$
Similarly, in maximal dimension $d=N-1$, the exponents
$w_{N-1,1}(\Theta)$ and $\hat{w}_{N-1,1}(\Theta)$ measure
the approximation to $\Theta$ by rational hyperplanes. Then, the
equalities
$$
w_{N-1,1}(\Theta) = \om_{n,m}(A) \and \hat{w}_{N-1,1}(\Theta) =
\omc_{n,m}(A)
$$
hold.

We address the following

\proclaim Problem 4.  Let $\delta$ and $N$ be integers with $0 \le
\delta \le  N-1$. Find the spectrum of the array of  exponents
$$
\Big\{ \dots,  w_{d,k}(\Theta) ,  \hat{w}_{d,k}(\Theta), \dots   ;
0\le d< N-1, 1 \le k \le  \min(\delta +1,N-\delta,d+1,N-d)  \Big\},
$$
where $\Theta$ ranges over the set of $\delta$-dimensional subvarieties
of $\bP^{N}(\bR)$.

For instance, when $\delta=0, N=3$ (approximating points in $\bP^3$),
we have to investigate the possible values of six exponents,
while for $\delta =1, N=3$ (approximating lines  in $\bP^3$), we have
to look at 8-tuples.  Notice that the spectra associated
to the pairs of dimensions $(\delta, N)$ and $(N-1 -\delta, N)$ are
permuted by duality.

Theorem 7 provides the answer of Problem 4
when $N=2$. Informations about the above array of exponents, including
lower bounds and transference results, may be deduced from
\cite{SchmB}. However, as far as we are aware, the generic value of
this array (if it does exist) remains unclear, unless when $\delta =0$
or when $\delta =N-1$.

\section{5. Approximation of dependent quantities}

In this section, we assume that either $m$ or $n$
is equal to 1. In other words, we let $\xi_1, \ldots , \xi_n$ be
$n$ real numbers, and we take for $A$ either the matrix
$(\xi_1, \ldots, \xi_n)$ or the matrix ${}^t (\xi_1, \ldots, \xi_n)$.
In Section 2, we have assumed that the $\xi_j$'s are independent.
We deal now with the more complicated situation of dependent $\xi_j$'s;
typically, we assume that $(\xi_1, \ldots, \xi_n)$ belong to some
given manifold. There is a broad literature on this subject, and we
direct the reader to the book of Bernik \& Dodson \cite{BeDo} for
results and many bibliographical references. In the sequel, we restrict
our attention to the case where $\xi_j = \xi^j$ for a given real number
$\xi$ and any integer $j$ with $1 \le j \le n$.

This situation is the most classical one. Indeed, in order to define in
1932
his classification of the real numbers $\xi$, Mahler \cite{Mah32}
introduced
the exponents of Diophantine approximation $w_n(\xi)$, which correspond
to the exponents $\om_{1,n}((\xi, \xi^2, \ldots, \xi^n))$
defined in Section 2, when $\xi$ is not algebraic of degree at most 
$n$.
In view of the particular questions
investigated in the present section, we do not keep the notation of
Definition 1, and we rather use the classical notation, recalled
in Definition 5 below.

\proclaim Definition 5.
Let $n \ge 1$ be an integer and let $\xi$ be a real number.
We denote by $w_n(\xi)$ (resp. by $\hw_n(\xi)$)
the supremum of the real numbers $w$ such that, for
arbitrarily large real numbers $X$
(resp. any sufficiently large real number $X$), the inequalities
$$
0 < |x_n \xi^n + \ldots + x_1 \xi + x_0| \le X^{-w}, \qquad  \max_{0
\le m \le n} \, |x_m| \le X,
$$
have a solution in integers $x_0, \ldots, x_n$.
We denote by $\la_n(\xi)$ (resp. by $ \hat{\lambda}_n(\xi)$)
the supremum  of the real numbers
$\lambda$ such that, for
arbitrarily large real numbers $X$
(resp. any sufficiently large real number $X$), the inequalities
$$
0 < |x_0| \le X, \qquad \max_{1 \le m \le n} \,
|x_0 \xi^m - x_m| \le X^{-\lambda},
$$
have a solution in integers $x_0, \ldots, x_n$.

Observe that $\la_n (\xi) = \om_{n, 1} ({}^t (\xi, \xi^2, \ldots,
\xi^n))$
and $\hat{\lambda}_n (\xi) = \omc_{n, 1} ({}^t (\xi, \xi^2, \ldots,
\xi^n))$
hold for any $n\ge 1$ and any real number $\xi$ not algebraic of degree
at most $n$.

Unfortunately, Theorem 2 and 7 do not imply any information regarding
the
values  of the functions $w_n$ and $\la_n$.
Solving a long-standing conjecture of Mahler,
Sprind\v zuk \cite{Spr69} proved in 1965 that $w_n (\xi) = n $
holds for any $n\ge 1$ for
almost all (with respect to the Lebesgue measure) real numbers $\xi$.
By (6), this implies that
$\la_n (\xi) = 1/n$ holds for any $n\ge 1$, for almost all real numbers
$\xi$.

Furthermore, it follows from the Schmidt Subspace Theorem
(see e.g. \cite{SchmLN}) that
$$
\hw_n (\xi) = w_n (\xi) = 1/ \la_n (\xi)
= 1/\hla_n (\xi) = \min\{d-1, n\}  \eqno (12)
$$
hold for any real algebraic number $\xi$ of degree $d$.
Thus, to investigate the
sets of values taken the functions $w_n$, $\hw_n$, $\la_n$ and
$\hla_n$, we need only to consider them on transcendental numbers.

\proclaim Definition 6. By spectrum of the function $w_n$
(resp. $\hw_n$, $\la_n$ and $\hla_n$), we mean the set of values
taken by $w_n$ (resp. $\hw_n$, $\la_n$ and $\hla_n$) on the
set of transcendental real numbers.

It is possible to construct explicit examples of real numbers $\xi$
with $w_n(\xi) = w$, for any given real number 
$w>(2n+1+\sqrt{4n^2 +1})/2$, see 
Theorem 7.7 \cite{BuLiv} 
for
references. As in the proof of Theorem 3, the theory of
Hausdorff dimension is a crucial tool for determining the spectrum
of $w_n$, a problem solved in 1983 by Bernik \cite{Ber}.

\proclaim Theorem 9. For any positive integer $n$, the spectrum
of $w_n$ is equal to $[n, + \infty]$.

\pro It follows from \cite{Ber} that, for any real number $\tau \ge n$,
we have
$$
\dim_{\rm H} \{\xi \in \R : w_n (\xi) = 
\tau (n+1) - 1\} = {1 \over \tau},
$$
where $\dim_{\rm H}$ denotes the Hausdorff dimension.  
We finish the proof of the theorem by observing that any
Liouville number $\xi$ satisfies $w_n (\xi)= w_1 (\xi) = + \infty$.
\cqfd

As far as we are aware, the spectra of the functions $\la_n$ have not
been studied up to now when $n \ge 2$.

\proclaim Theorem 10. For any integer $n \ge 1$, the spectrum  
of $\la_n$ includes  the interval   
$((3+\sqrt{4n^2 +1})/(2n),+\infty]$.   

\pro  
Notice that  the upper bound 
$$
\max_{1\le k \le n} | q^n\xi^k -q^{n-k}p^{k}|  \le 
n\max\{1,|\xi|\}^{n-1}\max\{| p| , |q|\}^{n-1}|q\xi-p|, 
$$
which holds for all integers $p$ and $q$,   implies the lower bound 
$$
\lambda_n(\xi) \ge { w_1(\xi) -n+1\over n}. 
$$
On the other hand, Khintchine's transference principle (6) provides us 
with the upper bound (see also Theorem 3.9 of \cite{BuLiv}) 
$$
\lambda_n(\xi) \le { w_n(\xi) -n +1\over n}. 
$$
Now, Theorem 7.7 in \cite {BuLiv} asserts that for any given real 
number $w > (2n+1+\sqrt{4n^2 +1})/2$, there exists a real number $\xi$ 
(that can be given explicitly) with  
$w_1(\xi) = w_n(\xi) = w$.  Then,  the equality 
$$
\lambda_n(\xi ) = {w - n  +1 \over n} 
$$
holds.  
\cqfd

\proclaim Problem 5. Let $n \ge 1$ be an integer. Is the   
spectrum of the function $\la_n$ equal to $[1/n,+\infty]$  ? 

We now turn our attention
to the exponents of uniform approximation $\hw_n$ and $\hla_n$,
introduced explicitly for the first time in \cite{BuLaA}, but already
studied
by Davenport \& Schmidt \cite{DaSc} in 1969.

\proclaim Theorem 11. For any integer $n \ge 1$ and any transcendental
real number $\xi$, we have
$$
\hla_n (\xi) \le 1 / \lceil n/2  \rceil
\quad {\it and} \quad \hw_n (\xi) \le 2n - 1.  \eqno (13)
$$
Furthermore, we have
$$
\hla_2 (\xi) \le (\sqrt{5} - 1)/2
\quad {\it and} \quad \hw_2 (\xi) \le (3 + \sqrt{5})/2.  \eqno (14)
$$

\pro The two bounds (13) combine results by
Davenport \& Schmidt \cite{DaSc} and  Laurent \cite{Lau}. We now sketch
a new proof of (14),
as a consequence of Theorems 6 and 7. We first establish that
$\lambda_2(\xi)\le 1$ whenever $\hla_2(\xi) > 1/2$, which   obviously 
may be assumed. 
Let $\epsilon >0$ and let
$\ux = (x_0,x_1,x_2)$ be a non-zero integer triple with large norm
$|\ux |=X$ such that
$$
\max\{|x_0 \xi -x_1|,|x_0 \xi^2 -x_2| \} \le
X^{-\lambda_2(\xi)+\epsilon}.
$$
It is shown in Lemma 2 of \cite{DaSc} that we may suppose without loss 
of   
generality that the Hankel determinant
$x_1^2- x_0x_2$ is non-zero. Arguing as Davenport \& Schmidt, we deduce
   the upper bound
$\lambda_2(\xi) \le 1+\epsilon$ from the estimate
$$
1 \le | x_1^2-x_0x_2| \le 2 X^{1-\lambda_2(\xi) +\epsilon},
$$
which is valid for arbitrarily large values of  $X$.
Using now  Corollary 2, we find the
inequalities
$$
1 \ge \lambda_2(\xi) \ge {\hla_2(\xi)^2\over 1 -\hla_2(\xi)},
$$
from which follows the first upper bound of (14). Notice finally that
the two upper bounds of (14) are equivalent
by Theorem 6. An alternative  proof of  the bound $\omc_2 (\xi) \le (3
+ \sqrt{5})/2$ may also be found in Arbour \& Roy \cite{ArRo}.
\cqfd

For a long time, it was believed that the upper bounds in (13)
could possibly be improved to $1/n$ and $n$, respectively.
This is, however, not true for $n=2$, as recently proved by Roy in a
series
of remarkable papers: the upper bounds given in (14) are best possible.
Section 6 is devoted to Roy's recent works and
some of their extensions.

\section{ 6. Some computations of exponents}

In this section, we restrict our attention to the functions
$\hw_2$ and $\hla_2$, which, as first shown by
Roy \cite{RoyA,RoyB}, take values
strictly larger than $2$ and $1/2$,
respectively, at some transcendental points.

\proclaim Theorem 12. There are real transcendental numbers $\xi$
with $\hw_2 (\xi) = (3 + \sqrt{5})/2$ and
$\hla_2 (\xi) = (\sqrt{5} - 1)/2$.

Recall that real numbers $\xi$ with either $\hw_n (\xi) > n$ or
$\hla_n (\xi) > 1/n$ for some $n \ge 2$ are transcendental, by 
Schmidt's
Subspace Theorem.

Actually, the result of Roy is slightly more precise, since he
constructed
numbers $\xi$ for which there exists some positive constant $c$ such
that
the systems of inequalities
$$
\eqalign{
|x_2 \xi^2 + x_1 \xi + x_0| & \le c \, X^{-(3 + \sqrt{5})/2}, \cr
|x_1|, |x_2| & \le X, \cr} \eqno (15)
$$
and
$$
\eqalign{
|x'_0 \xi + x'_1| & \le c \, X^{-(\sqrt{5} - 1 )/2}, \cr
|x'_0 \xi^2 + x'_2| & \le c \, X^{-(\sqrt{5} - 1 )/2}, \cr
|x'_0| & \le X, \cr}  \eqno (16)
$$
have non-zero integer solutions $(x_0, x_1, x_2)$ and
$(x'_0, x'_1, x'_2)$, respectively, for any real
number $X>1$. Such a result is quite surprising, since the volumes of
the convex bodies defined by (15) and (16) tend rapidly to zero as $X$
grows to infinity. According to Roy, such a real
number $\xi$ is called an extremal number. The set
of extremal numbers, which is countable \cite{RoyB}, has been further
studied in \cite{RoyF}.

\medskip

Let us now give a real number with these extremal properties.
Let $\{a, b\}^*$ denote the monoid of words on
the alphabet $\{a, b\}$ for the product given by the
concatenation. The Fibonacci sequence in $\{a, b\}^*$ is the sequence
of words $(f_i)_{i\ge 0}$ defined recursively by
$$
f_0 = b, \ens f_1 = a, \quad {\rm and} \quad
f_i = f_{i-1} f_{i-2} \, (i \ge 2).
$$
Since, for every $i \ge 1$, the word $f_i$ is a prefix of $f_{i+1}$,
this sequence converges to an infinite word $f = abaabab\ldots $ called
the Fibonacci word
on $\{a, b\}$. For two positive distinct integers
$a$ and $b$, let $\xi_{a, b} = [0; a, b, a, a, b, a, \ldots]$ be the
real number whose sequence of partial quotients is given by the letters
of the Fibonacci word on $\{a, b\}$. Then, $\xi_{a, b}$ satisfies the
properties stated in Theorem 12.

\medskip

\noi {\bf Sketch of the proof of Theorem 12.}

Following \cite{RoyA}, we show that there exists a suitable
constant $c$ for which the system (16) with $\xi = \xi_{a, b}$ has
a non-zero integer solution for any real number $X >1$.

We begin with a property of the Fibonacci word $f$.
Let $(F_m)_{m \ge 0}$ denote the Fibonacci sequence defined by
$F_0 = 0$, $F_1 = 1$, and $F_{m + 2} = F_{m+1} + F_m$, for $m\ge 0$.
It is well known (see e.g. \cite{ADQZ}) that the sequence $(\phi_n)_{n
\ge 2}$
formed by the prefixes $\phi_n$
of $f$ of length $ F_{n+2} - 2$
has the following property: for any $n\ge 2$, the word $\phi_n$
is a palindrome. Observe that $\phi_2 = a$, $\phi_3 = aba$,
and $\phi_4 = abaaba$. Furthermore, we have
$$
\phi_n = \phi_{n-1} ab \phi_{n-2}, \qquad \hbox{for $n \ge 4$ even},
\eqno (17)
$$
and
$$
\phi_n = \phi_{n-1} ba \phi_{n-2}, \qquad \hbox{for $n \ge 5$ odd}.
\eqno (18)
$$

Before going on with the proof, we make the following observation,
extracted from \cite{AdBu}.
Let $\eta = [0; a_1, a_2, \ldots]$ be a positive real irrational  
number, and
denote by $p_n/q_n$ its convergents, that is,
$p_n/q_n = [0; a_1, a_2, \ldots, a_{n}]$. By the theory of  
continued fraction, we have
$$
M_n := \pmatrix{ q_n & q_{n-1} \cr p_n & p_{n-1} \cr} =
\pmatrix{ a_1 & 1 \cr 1 & 0 \cr}
\pmatrix{ a_2 & 1 \cr 1 & 0 \cr} \ldots
\pmatrix{ a_n & 1 \cr 1 & 0 \cr},  
$$
and, since such a decomposition is unique,
the matrix $M_n$ is symmetrical if, and only if,
the word $a_1 a_2 \ldots a_n$ is a palindrome, that is, if, and
only if,
we have $a_j = a_{n+1-j}$ for any integer $j$ with $1 \le j \le n$.  
In this case, we have $p_n = q_{n-1}$ and, by the
theory of continued fraction, we get
$$
\biggl| \eta - {p_n \over q_n} \biggr| < {1 \over q_n^2}
\quad {\rm and} \quad
\biggl| \eta - {p_{n-1} \over q_{n-1}} \biggr| < {1 \over q_{n-1}^2}.
$$
We then infer from $0 < \eta < 1$, $a_1 = a_n$ and  
$q_n \le (a_n + 1) q_{n-1}$ that
$$
\eqalign{
\biggl| \eta^2 - {p_{n-1} \over q_n} \biggr| & \le
\biggl| \eta^2 - {p_{n-1} \over q_{n-1}} \cdot  {p_n \over q_n} \biggr|
\le \biggl| \eta  + {p_{n-1} \over q_{n-1}} \biggr|
\cdot \biggl| \eta -  {p_n \over q_n}  \biggr|
+ {1 \over q_n q_{n-1}} \cr
& \le 2 \biggl| \eta  - {p_n \over q_n} \biggr| + {1 \over q_n q_{n-1}}
< {a_1 + 3 \over q_n^2}. \cr}  
$$
Consequently, if the sequence of the partial quotients of $\eta$
is bounded and begins with
infinitely many palindromes, then $\eta$ and $\eta^2$ are
simultaneously very well approximable by rational
numbers of the same denominator, we have $\la_2 (\eta) = 1$ and $\eta$
is either quadratic, or transcendental, by (12).
As noted in
\cite{AdBu,AdBuAMM}, this observation gives a very short proof of the
transcendence of the Thue--Morse continued fraction, first established
by Queff\'elec \cite{Que} (see also \cite{AlSh}).

\medskip

For any $n\ge 2$, denote by $Q_n$ the
denominator of the rational number whose partial quotients are given by
the letters of $\phi_n$. The above observation shows that, for
a suitable constant $c_1$ and any $n\ge 4$, the system
$$
\eqalign{
|x_0 \xi_{a, b} + x_1| & \le c_1 \, Q_n^{-1}, \cr
|x_0 \xi_{a, b}^2 + x_2| & \le c_1 \, Q_n^{-1}, \cr
|x_0| & \le Q_n, \cr}
$$
has a non-zero solution, that we denote by
$(x_0^{(n)}, x_1^{(n)}, x_2^{(n)})$. Observe that $Q_n = x_0^{(n)}$.

Furthermore, it follows from (17) and (18) that
$$
\pmatrix{ x_0^{(n)} & x_1^{(n)} \cr x_1^{(n)} & x_2^{(n)} \cr} =
\pmatrix{ x_0^{(n-1)} & x_1^{(n-1)} \cr x_1^{(n-1)} & x_2^{(n-1)} \cr}
\times S_n \times
\pmatrix{ x_0^{(n-2)} & x_1^{(n-2)} \cr x_1^{(n-2)} & x_2^{(n-2)} \cr},
$$
where
$$
S_n = \pmatrix{a & 1 \cr 1 & 0 \cr} \times \pmatrix{b & 1 \cr 1 & 0 
\cr}
\qquad {\rm or} \qquad \pmatrix{b & 1 \cr 1 & 0 \cr}
\times \pmatrix{a & 1 \cr 1 & 0 \cr},
$$
according as $n$ is even or odd. This yields
$$
x_0^{(n)} = \pmatrix{x_0^{(n-1)} & x_1^{(n-1)} \cr} \times S_n \times
\pmatrix{ x_0^{(n-2)} \cr x_1^{(n-2)}  \cr},
$$
which implies
$$
\eqalign{
\lim_{n \to + \infty} \, {Q_n \over Q_{n-1} Q_{n-2}} & =
\pmatrix{1 & \xi_{a,b} \cr} \times \pmatrix{a & 1 \cr 1 & 0 \cr} \times
\pmatrix{b & 1 \cr 1 & 0 \cr} \times \pmatrix{1 \cr \xi_{a,b} \cr} \cr
& = \xi_{a,b}^2 + (a+b) \xi_{a,b} + (ab+1). \cr}  \eqno (19)
$$
Set $\gamma = (1 + \sqrt{5})/2$ and
$\kappa_n = Q_n \, Q_{n-1}^{-\gamma}$, for any integer $n\ge 2$.
Since $\gamma = 1 + 1/\gamma$, we get
$$
\kappa_n = {Q_n \over Q_{n-1} Q_{n-2}} \, \kappa_{n-1}^{-1/\gamma},
$$
thus, by (19), there exist positive constants $c_3 > c_2$ such that
$$
c_2 \kappa_{n-1}^{-1/\gamma} < \kappa_n < c_3 \kappa_{n-1}^{-1/\gamma},
\qquad \hbox{for any $n\ge 2$}.
$$
By induction, this yields
$$
c_4 \, Q_n^{(1+\sqrt{5})/2} \le Q_{n+1} \le c_5 Q_n^{(1+\sqrt{5})/2},
\qquad \hbox{for any $n \ge 4$}, \eqno (20)
$$
with $c_4 = \min\{Q_2, c_2^{\gamma}/c_3\}$
and $c_5 = \max\{Q_2, c_3^{\gamma}/c_2\}$.

\medskip

Let then $X$ be a sufficiently large real number.
There exists an integer $n \ge 4$ such that  
$Q_n \le X < Q_{n+1}$. The system
$$
\eqalign{
|x_0 \xi_{a, b} + x_1| & \le X^{-1/2}, \cr
|x_0 \xi_{a, b}^2 + x_2| & \le X^{-1/2}, \cr
|x_0| & \le X \cr}
$$
has the non-zero integer solution $(x_0^{(n)}, x_1^{(n)}, x_2^{(n)})$,
and we have
$$
\max\{|x_0^{(n)} \xi_{a, b} + x_1^{(n)}|, |x_0^{(n)} \xi_{a, b}^2 +
x_2^{(n)}|\}
\le c_1 \, Q_n^{-1} \le c_1 \, X^{- (\log Q_n)/(\log X)}
\le c_6 \, X^{-2/(1+\sqrt{5})},
$$
for a suitable positive constant $c_6$, by (20).
This shows that $\hla_2(\xi_{a, b}) \ge (\sqrt{5}-1)/2$.
By Theorems 6 and 11, this yields Theorem 12. \cqfd   

A natural question is now the study of the spectra of the functions
$\hw_2$ and $\hla_2$. Roy \cite{RoyB}
proved that there are only countably many
real numbers $\xi$ for which the system (15) (resp. (16)) has a 
non-zero
solution for any $X>1$. Shortly thereafter,
Bugeaud \& Laurent \cite{BuLaA} extended
Roy's construction and found uncountably many values
taken by $\hw_2$ and $\hla_2$.

\proclaim Theorem 13. Let $(s_j)_{j \ge 1}$ be a bounded sequence of
integers
and set
$$
\sigma = \liminf_{k} \, [0; s_k,s_{k-1}, \dots, s_1].
$$
There exist real numbers $\xi$ with
$$
\displaylines{
\lambda_2 (\xi) = {1 }, \qquad
w_2(\xi)  = 1 + {2 \over \sigma},
\cr
\hat{\lambda}_2(\xi)= {1 +\sigma \over 2 +\sigma}, \qquad
\hw_2(\xi) = 2 +\sigma.
\cr}
$$
In particular, the spectra of $\hw_2$ and $\hla_2$ have Hausdorff
dimension 1.

Notice that these four exponents satisfy the relation
$$
\lambda_2(\xi) = { w_2(\xi)(\hat{w}_2(\xi) -1)\over w_2(\xi)+
\hat{w}_2(\xi)}.
$$
Thus,  the numbers $\xi$ constructed in   Theorem 13,   provide us with
   examples of extremal matrices $A= (\xi,\xi^2)$ for which
the  lower bound $v' \ge v(w-1)/(v+w)$ given by  Theorem 7  turns out
to be an equality.
   Theorem 13 is proved in \cite{BuLaA}, and we refer to Cassaigne
\cite{Cass}
and to \cite{BuLaA} for further results on the function
$$
(s_j)_{j \ge 1} \longmapsto \liminf_{k} \, [0; s_k,s_{k-1}, \dots, 
s_1].
$$

It is tempting to believe that the spectra of $\hw_2$ and $\hla_2$ 
enjoy
a structure of `Markoff spectrum', and that  
$(3+\sqrt{5})/2$ and $(\sqrt{5} - 1)/2$ are
isolated points of the spectra of $\hw_2$ and $\hla_2$, respectively.
This is, however, not true, as recently established by Roy \cite{RoyE}.

\proclaim Theorem 14. The spectra of $\hw_2$ and $\hla_2$ are dense
in $[2, (3+\sqrt{5})/2]$ and $[1/2, (\sqrt{5} - 1)/2]$,
respectively.

To prove Theorem 14, Roy produces countably many
real numbers $\xi$ of `Fibonacci
type' by suitably modifying the constructions of
Section 6 of \cite{RoyB} and Section 5 of \cite{RoyD}. In view of
results from \cite{BuLaA},
one may ask whether there exist transcendental real numbers $\xi$
not of that type which satisfy $\hat{w}_2(\xi) > 1+\sqrt{2}$.

\proclaim Problem 6. Determine the spectra of
the functions $\hw_2$ and $\hla_2$.

%equal to $[2, (3+\sqrt{5})/2]$ and $[1/2, (\sqrt{5} - 1)/2]$,
%respectively?

Works of Fischler \cite{FiscA,FiscB,FiscC} bring some
light on Problem 6. Furthermore, he informed us that he established
the existence of a (small) explicitly computable positive number
$\epsilon$ such that the intersection of the spectrum of $\hw_2$
with $[(3+\sqrt{5})/2 - \epsilon, (3+\sqrt{5})/2]$ is countable.
Consequently, the spectrum of $\hw_2$ is not equal to the whole
interval $[2, (3+\sqrt{5})/2]$.

\vskip 28mm

\centerline{\bf References}

\vskip 5mm

\beginthebibliography{999}

\bibitem{AdBu}
B. Adamczewski and Y. Bugeaud,
{\it Palindromic continued fractions}.
Preprint.

\bibitem{AdBuAMM}
B. Adamczewski and Y. Bugeaud,
{\it A short proof of the
transcendence of the Thue--Morse continued fractions},
Amer. Math. Monthly. To appear.

\bibitem{ADQZ}
J.-P. Allouche, J. L. Davison, M. Queff\'elec, and L. Q. Zamboni,
{\it Transcendence of Sturmian or morphic continued fractions},
J. Number Theory 91 (2001), 39--66.

\bibitem{AlSh}
J.-P. Allouche and J. Shallit,   
Automatic Sequences: Theory, Applications, Generalizations,
Cambridge University Press, Cambridge, 2003.

\bibitem{Apf}
A. Apfelbeck,
{\it A contribution to Khintchine's principle of transfer},
Czechoslovak Math. J. 1 (1951), 119--147.

\bibitem{ArRo}
B. Arbour and D. Roy,
{\it A Gel'fond type criterion in degree two},
Acta Arith. 111 (2004), 97--103.

\bibitem{BDV}
V. V. Beresnevich, H. Dickinson, and S. L. Velani,
{\it Sets of exact `logarithmic order' in
the theory of Diophantine approximation}, Math. Ann. 321 (2001),
253--273.

\bibitem{Ber}
V. I. Bernik,
{\it Application of the Hausdorff dimension in the theory of
Diophantine approximations}, Acta Arith. 42 (1983), 219--253 (in
Russian).
English transl. in Amer. Math. Soc. Transl. 140 (1988), 15--44.

\bibitem{BeDo}
V. I. Bernik and M. M. Dodson,
Metric Diophantine approximation on
manifolds, Cambridge Tracts in Mathematics 137, Cambridge University
Press,
1999.

\bibitem{BuLiv}
Y. Bugeaud,
Approximation by algebraic numbers,
Cambridge Tracts in Mathematics, Cambridge, 2004.

\bibitem{BuLaA}
Y. Bugeaud and M. Laurent,
{\it Exponents of Diophantine Approximation and
Sturmian Continued Fractions},
Ann. Inst. Fourier (Grenoble) 55 (2005), 773--804.

\bibitem{BuLaB}
Y. Bugeaud and M. Laurent,
{\it Exponents of homogeneous and inhomogeneous Diophantine
Approximation},
Moscow Math. J. To appear.

\bibitem{Cass}
J. Cassaigne,
{\it Limit values of the recurrence quotient of Sturmian sequences},
Theor. Comput. Sci. 218 (1999), 3--12.

\bibitem{Cas}
J. W. S. Cassels,
An introduction to Diophantine Approximation,
Cambridge Tracts in Math. and Math. Phys., vol. 99, Cambridge
University Press, 1957.

\bibitem{ChLu}
C. Chabauty et \'E. Lutz,
{\it Sur les approximations diophantiennes lin\'eaires r\'eelles (I).
Probl\`eme homog\`ene},
C. R. Acad. Sci. Paris 231 (1950), 887--888.

\bibitem{DaSc}
        {H. Davenport and W. M. Schmidt},
{\it Approximation to real numbers by
algebraic integers}, Acta Arith. {15} (1969), 393--416.

\bibitem{DSa}
H. Davenport and W. M. Schmidt,
{\it Dirichlet's theorem on
Diophantine approximation}, Symposia Mathematica, Vol. IV
(INDAM, Rome, 1968/69), pp. 113--132, Academic Press, London, 1970.

\bibitem{DaScB}
H. Davenport and W. M. Schmidt,
{\it Dirichlet's theorem on
Diophantine approximation. II}, Acta Arith. {16} (1970), 413--423.

\bibitem{DV}
H. Dickinson and S. L. Velani,
{\it Hausdorff measure and
linear forms}, J. reine angew. Math. 490 (1997), 1--36.

\bibitem{Dy47}
F. J. Dyson,
{\it On simultaneous Diophantine approximations},
Proc. London Math. Soc. 49 (1947), 409--420.

\bibitem{FiscA}
S. Fischler,
{\it Spectres pour l'approximation d'un nombre r\'eel et de son
carr\'e},
C. R. Acad. Sci. Paris 339 (2004), 679--682.

\bibitem{FiscB}
S. Fischler,
{\it Palindromic Prefixes and Episturmian Words},
J. Combin. Theory, Series A. To appear.

\bibitem{FiscC}
S. Fischler,
{\it Palindromic Prefixes and Diophantine Approximation},
Monatsh. Math. To appear. 

\bibitem{GrLe}
P. M. Gruber and C. G. Lekkerkerker,
Geometry of numbers,
Series Bibliotheca Mathematica 8,
North--Holland, Amsterdam, 1987.

\bibitem{Jar35}
V. Jarn\'\i k,
{\it \"Uber einen Satz von A. Khintchine},
Pr\'ace Mat.-Fiz. 43 (1935), 1--16.

\bibitem{Jar35a}
V. Jarn\'\i k,
{\it O simult\'ann\i\' ch diofantick\'ych approximac\i\' ch},
Rozpravy T\'r. \v Cesk\'e Akad 45, c. 19 (1936), 16 p.

\bibitem{Jar36}
V. Jarn\'\i k,
{\it \"Uber einen Satz von A. Khintchine, 2. Mitteilung},
Acta Arith. 2 (1936), 1--22.

\bibitem{Jar38}
V. Jarn\'\i k,
{\it Zum Khintchineschen ``\"Ubertragungssatz''},
Trav. Inst. Math. Tbilissi 3 (1938), 193--212.

\bibitem{Jar50}
V. Jarn\'\i k,
{\it Une remarque sur les approximations diophantiennes lin\'eaires},
Acta Sci. Math. Szeged 12 (1950), 82--86.

\bibitem{Jar54}
V. Jarn\'\i k,
{\it Contribution \`a la th\'eorie des approximations diophantiennes
lin\'eaires et homog\`enes},
Czechoslovak Math. J. 4 (1954), 330--353 (in Russian, French summary).

\bibitem{Jar59}
V. Jarn\'\i k,
{\it Eine Bemerkung zum \"Ubertragungssatz},
B\u ulgar. Akad. Nauk Izv. Mat. Inst. 3 (1959), 169--175.

\bibitem{Jar59b}
V. Jarn\'\i k,
{\it Eine Bemerkung \"uber diophantische Approximationen},
Math. Z. 72 (1959), 187--191.

\bibitem{Kh25}
A. Ya. Khintchine,
{\it Zwei Bemerkungen zu einer Arbeit des Herrn Perron},
Math. Z. 22 (1925), 274--284.

\bibitem{Kh26b}
A. Ya. Khintchine,
{\it \"Uber eine Klasse linearer diophantischer Approximationen},
Rendiconti Circ. Mat. Palermo 50 (1926), 170--195.

\bibitem{KhB}
A. Ya. Khintchine,
{\it Regular systems of linear equations and a general problem of
\v Ceby\v sev}, Izvestiya Akad. Nauk SSSR Ser. Mat. 12 (1948),
249--258.

\bibitem{KlWe}
D. Kleinbock and B. Weiss,   
{\it Dirichlet's theorem on Diophantine approximation and homogeneous flows}.
Preprint.

\bibitem{Lau}
        M. Laurent,
{\it Simultaneous rational approximation to the successive powers
of a real number},
Indag. Math. 11 (2003), 45--53.

\bibitem{LauB}
        M. Laurent,
{\it Exponents of Diophantine Approximation in dimension two},
Canad. J. Math. To appear.

\bibitem{Les}
J. Lesca,
{\it Sur un r\'esultat de Jarn\'\i k},
Acta Arith. 11 (1966), 359--364.

\bibitem{Mah32}
K. Mahler,
{\it Zur Approximation der Exponentialfunktionen und des
Logarithmus. I, II},
J. reine angew. Math. 166 (1932), 118--150.

\bibitem{Que}
M. Queff\'elec,
{\it Transcendance des fractions continues de Thue--Morse},
J. Number Theory 73 (1998), 201--211.

\bibitem{Rand}
H. Randriambololona,
{\it Hauteurs des sous-sch\'emas de dimension nulle},
Ann. Inst. Fourier 53 (2003), 2155--2224.

\bibitem{RoyA}
D. Roy,
{\it Approximation simultan\'ee d'un nombre et son carr\'e},
C. R. Acad. Sci. Paris 336 (2003), 1--6.

\bibitem{RoyB}
        D. Roy,
{\it Approximation to real numbers by cubic algebraic numbers, I},
Proc. London Math. Soc. 88 (2004), 42--62.

\bibitem{RoyC}
        D. Roy,
{\it Approximation to real numbers by cubic algebraic numbers, II},
Annals of Math. 158 (2003), 1081--1087.

\bibitem{RoyD}
D. Roy,
{\it Diophantine approximation in small degree},
Number Theory, 269--285, CRM Proc. Lecture Notes 36,
Amer. Math. Soc., Providence, RI, 2004.

\bibitem{RoyE}
D. Roy,
{\it On two exponents of approximation related to a real number and
its square},
Canad. J. Math. To appear.

\bibitem{RoyF}
D. Roy,
{\it On the continued fraction expansion of a class of numbers},  
Proceedings of the Conference in honour of the 70th birthday 
of W. M. Schmidt, European Mathematical Society. To appear.

\bibitem{SchmLN}
W. M. Schmidt,
Diophantine Approximation, Lecture Notes in Math.
{785}, Springer, Berlin, 1980.

\bibitem{SchmB}
W. M. Schmidt,
{ \it On heights of algebraic subspaces and diophantine 
approximations},
Annals of Math. 85 (1967), 430--472.

\bibitem{Spr69}
{V. G. Sprind\v zuk},
Mahler's problem in metric number theory,
Izdat. ``Nauka i Tehnika'', Minsk, 1967 (in Russian).
English translation by B. Volkmann, Translations of Mathematical
Monographs,
Vol. 25, American Mathematical Society, Providence, R.I., 1969.

\endthebibliography

\vskip1cm

\noindent Yann Bugeaud  \hfill{Michel Laurent}

\noindent Universit\'e Louis Pasteur
\hfill{Institut de Math\'ematiques de Luminy}

\noindent U. F. R. de math\'ematiques
\hfill{C.N.R.S. -  U.M.R. 6206 - case 907}

\noindent 7, rue Ren\'e Descartes      \hfill{163, avenue de Luminy}

\noindent 67084 STRASBOURG  (FRANCE)
\hfill{13288 MARSEILLE CEDEX 9  (FRANCE)}

\vskip2mm

\noindent {\tt bugeaud@math.u-strasbg.fr}
\hfill{{\tt laurent@iml.univ-mrs.fr}}

\bye